\documentclass[11pt]{article}
\usepackage{amsmath,amssymb,amsthm}
\usepackage{booktabs,graphicx,tikz,geometry,placeins}
\usepackage{float}
\usepackage[hidelinks]{hyperref}
\newtheorem{theorem}{Theorem}
\newtheorem{proposition}[theorem]{Proposition}

\theoremstyle{remark}

\title{Why $2$, $3$, and $5$? Lambert Series, Prime Logarithms, and a $48$-Step Octave for Musical Tuning}
\author{Shubho Banerjee\\
\small Professor Emeritus of Physics\\
\small Rhodes College, Memphis, Tennessee, USA\\
\small \texttt{banerjees@rhodes.edu}}
\date{}

\begin{document}
\maketitle

\begin{abstract}
The primes $2$, $3$, and $5$ generate many familiar
frequency ratios in musical tuning. We identify
three mathematical properties that distinguish
these primes in a construction based on modular
transformations of Lambert series. We derive an
exact identity for $\log p$ for every odd prime
$p$, together with a separate identity for
$\log2$. Within this construction, all transformed
Lambert-series arguments can be chosen real
precisely for $p=2,3,5$. In the identities
obtained for these three primes, every series
argument has magnitude at most $e^{-2\pi}$,
whereas the identity for $\log7$ requires
complex arguments and contains a more slowly
converging series.

The leading $\pi$-terms for $\log2$, $\log3$, and
$\log5$ generate a $48$-step octave grid and
recover its conventional prime mapping
$(48,76,111)$ exactly, without rounding, while
including the leading $\pi$-term for $\log7$
enlarges the grid to $1680$ steps. The identity
for $\log7$ introduces no auxiliary logarithms
of larger primes, whereas for every prime
$11\le p<20000$, an exact computation
identifies such logarithms in the construction.

The Lambert-series remainders provide exact
corrections for just intervals and musical
commas. For a published theoretical model of
$22$ Hindustani shrutis, the additive
$48$-step mapping assigns all $22$ ratios
to distinct positions, with corrections no
greater than $15.3$ cents in magnitude.
These results establish precise arithmetic
distinctions that parallel the prominence
of $5$-limit tuning, without implying that
number theory determines musical practice.
\end{abstract}

\section{Introduction: the musical question}

A musical interval is a ratio of frequencies. Doubling a frequency gives an octave; the ratios $3/2$ and $5/4$ give a pure fifth and a pure major third. In \emph{$5$-limit just intonation}, musical intervals are constructed by multiplying powers of the primes $2$, $3$, and $5$. Such ratios can be written as
\begin{equation}\label{eq:five_limit}
 R=2^a3^b5^c.
\end{equation}
Here $a$, $b$, and $c$ are integers. The term \emph{limit} refers to the largest prime allowed in a ratio, not the number of notes in a scale. This framework is familiar in Western tuning theory \cite[Ch.~5]{Benson2006} (see also \cite{Scholtz1998}); related simple-ratio models also appear in discussions of Hindustani music and its traditional \emph{shrutis} \cite{Thakur2015}. Measured intonation in living musical traditions is more flexible than any fixed list of fractions \cite{JairazbhoyStone1963,SerraEtAl2011}.

The prime $7$ also gives musically meaningful intervals, such as the harmonic seventh $7/4$. It is not absent from musical history \cite[Ch.~6]{Benson2006}. Nevertheless, the widespread use of tuning constructions generated by $2$, $3$, and $5$ motivates a mathematical question: does this particular set of primes have a distinguishing property in a number-theoretic construction? We explore one such property without suggesting that it caused the historical development of any tuning system.

Our starting point is a family of identities expressing $\log p$ as a rational multiple of $\pi$ together with rapidly converging Lambert series. Two features make the identities useful here. First, the method treats prime numbers systematically: it produces a formula for every odd prime, so comparisons with $7$, $11$, and $13$ are not arbitrary. Second, because musical pitch in cents is proportional to a logarithm, the same identities separate a simple tempered pitch from an exact correction.

Two familiar discrepancies in musical tuning illustrate the role of the exact corrections. Twelve pure fifths do not quite equal seven octaves; their difference is the \emph{Pythagorean comma} $3^{12}/2^{19}$. Similarly, the Pythagorean major third $81/64$ differs from the pure major third $5/4$ by the \emph{syntonic comma} $81/80$. In our construction, the Pythagorean comma arises entirely from Lambert-series corrections, whereas the syntonic comma consists of one step on the $48$-step grid together with a Lambert-series correction. Both discrepancies are therefore accounted for exactly.

The exposition is organized for two kinds of reading. Sections~2 and~3 state the arithmetic results and explain their significance without requiring the reader to verify each step. Section~4 applies them to intervals, cents, and musical commas. Full derivations and a reproducible computational test appear in Appendix~A.

\section{From a Lambert series to a prime logarithm}

\subsection{The two identities behind the construction}

Define
\begin{equation}\label{eq:L}
 \mathcal L(q)=\sum_{n=1}^{\infty}\frac{q^n}{n(1-q^n)},\qquad |q|<1.
\end{equation}
For small $|q|$ the terms decay rapidly, making this series useful for high-precision identities. It is closely related to the Dedekind eta function through $\log\eta(\tau)=\pi i\tau/12-\mathcal L(e^{2\pi i\tau})$.

The construction uses two standard transformations, discussed in our earlier work \cite{BW2018}. For $\operatorname{Re}(t)>0$, \emph{modular inversion} gives
\begin{equation}\label{eq:inversion}
 \mathcal L(e^{-2\pi t})-\mathcal L(e^{-2\pi/t})
 =\frac12\log t-\frac{\pi}{12}\left(t-\frac1t\right).
\end{equation}
For a prime $p$, \emph{prime splitting} gives
\begin{equation}\label{eq:splitting}
 \sum_{j=0}^{p-1}\mathcal L\!\left(q^{1/p}e^{2\pi i j/p}\right)
 =(p+1)\mathcal L(q)-\mathcal L(q^p).
\end{equation}
The first replaces a series argument by its reciprocal modular parameter; the second combines the $p$ roots of a series argument. Neither identity refers to music. Their relevance appears when we use them to represent logarithms of the primes that generate musical ratios.

\subsection{One exact formula for every odd prime}

For an odd prime $p$, write $r=(p-1)/2$ and define
\begin{equation}\label{eq:Sp}
 S_p=\sum_{n=-r}^{r}\frac1{1+n^2},\qquad
 z_{p,n}=\exp\!\left(-\frac{2\pi p}{1+in}\right).
\end{equation}
The numbers $z_{p,n}$ all have absolute value less than one. They occur in conjugate pairs, so their total Lambert-series contribution is real even when individual arguments are complex.

\begin{theorem}[General prime-logarithm identity]\label{thm:prime}
For every odd prime $p$,
\begin{align}\label{eq:general}
\log p={}&\frac{\pi}{6}\left(S_p-\frac1p\right)
 +\frac1p\sum_{n=1}^{r}\log(1+n^2)
 -\frac{2(p+1)}p\mathcal L(e^{-2\pi})\nonumber\\
 &+\frac2p\mathcal L(e^{-2\pi p})
 +\frac2p\sum_{n=-r}^{r}\mathcal L(z_{p,n}).
\end{align}
\end{theorem}

What should a reader notice in this equation? It contains three types of terms. The rational multiple of $\pi$ provides an elementary component; the logarithms $\log(1+n^2)$ reveal which additional prime numbers enter the construction; and the Lambert-series terms complete the exact identity. When the auxiliary logarithms can be eliminated, the result consists of a rational multiple of $\pi$ and Lambert-series corrections. In the construction below, the $\pi$-terms provide useful approximations for $2$, $3$, and $5$, while the Lambert series restore exact equality. The identity remains valid for larger primes, including $11$ and $13$, even when their auxiliary logarithms cannot be eliminated recursively. Appendix~A contains the proof.

\subsection{The first three primes: three exact formulas}

The prime-logarithm identity provides one systematic construction for odd primes. The even prime $2$ is handled separately. The following three exact identities were obtained in our earlier work, Proposition~6.4 of \cite{BW2018}%
\footnote{In Proposition~6.4 of \cite{BW2018}, the coefficient of $\mathcal L(e^{-6\pi})$ in the formula for $\log3$ was printed as $16/3$. The correct coefficient is $4/3$, as used here.}:
\begin{align}
\log2={}&\frac{2\pi}{9}-\frac83\mathcal L(e^{-2\pi})
 +\frac83\mathcal L(e^{-4\pi}),\label{eq:log2}\\[3pt]
\log3={}&\frac{19\pi}{54}-\frac{32}{9}\mathcal L(e^{-2\pi})
 +\frac43\mathcal L(-e^{-3\pi})\nonumber\\
 &+\frac89\mathcal L(e^{-4\pi})
 +\frac43\mathcal L(e^{-6\pi}),\label{eq:log3}\\[3pt]
\log5={}&\frac{37\pi}{72}-\frac83\mathcal L(e^{-2\pi})
 +\frac23\mathcal L(e^{-4\pi})\nonumber\\
 &+\mathcal L(-e^{-5\pi})
 +\mathcal L(e^{-10\pi}).\label{eq:log5}
\end{align}

The equations have exactly the same structure: a rational multiple of $\pi$
plus a rapidly converging remainder. The remainder makes the formula exact;
the leading $\pi$-term provides a simple approximation. When these
three leading terms are measured in octaves, they give a $48$-step tempered
grid containing the usual $12$-tone equal-tempered positions. Their exact
Lambert-series remainders recover just-intonation ratios. Section~4
explains both facts in musical terms.

Identities \eqref{eq:log3} and \eqref{eq:log5} follow by applying Theorem~\ref{thm:prime} at $p=3$ and $p=5$, respectively, and eliminating the auxiliary $\log2$ terms with \eqref{eq:log2}. The same procedure produces the identity for $p=7$ in Section~3. The starting identity \eqref{eq:log2} is obtained by combining modular inversions to eliminate a slowly converging Lambert series. Other exact logarithmic identities exist, but once this starting identity and the symmetric representatives are specified, the identities obtained by the stated algorithm are determined algebraically. In particular, the rational multiples of $\pi$ in \eqref{eq:log3} and \eqref{eq:log5}, and hence the $48$-step grid examined below, arise from the stated procedure rather than from a separate selection of musical intervals.

\subsection{Why $2$, $3$, and $5$ are the real-argument cases}

The word \emph{real} here refers to the \emph{arguments of the Lambert series}, not to whether their sum or the logarithm is real. The imaginary part of $1/(1+in)$ causes a complex phase. If that phase is an integer multiple of $\pi$, the argument is real (possibly negative). This is exactly what happens for all the required arguments at $p=2,3,5$, but not for larger primes.

\begin{theorem}[The real-argument threshold]\label{thm:real}
Choose one integer representative $n$ of each residue class modulo a prime $p$ and form $z_{p,n}=\exp(-2\pi p/(1+in))$. All $p$ arguments can be made real if and only if $p\in\{2,3,5\}$.
\end{theorem}

For a reader interested in the main idea rather than its proof, the criterion is particularly simple:
\begin{equation}\label{eq:realcriterion}
 z_{p,n}\in\mathbb R\quad\Longleftrightarrow\quad
 (1+n^2)\mid 2p.
\end{equation}
For $p=5$, choose $n=-2,-1,0,1,2$: the numbers $1+n^2$ are $5,2,1,2,5$, and each divides $10$. At $p=7$, the residue class represented by $n=2$ already gives $1+n^2=5$, which does not divide $14$. Appendix~A proves that changing representatives cannot remove this obstruction for any prime $p\ge7$.

The theorem therefore singles out the same first three primes used to define the $5$-limit. It does \emph{not} show that primes above $5$ are musically unusable.

\section{What changes at $7$, $11$, and beyond?}

\subsection{A second threshold: whether the identity closes}

There is another distinction, separate from the reality of the arguments. To use \eqref{eq:general} as a recursive formula for $\log p$, we must deal with the logarithms of the integers $1+n^2$. Factoring these integers shows which prime logarithms appear. If they involve only smaller primes, or $p$ itself with a coefficient that can be moved to the left, the formula can be reduced using identities already obtained. We call this \emph{closure} of the symmetric construction.

For $p=7$ we need $n=1,2,3$, which give $2,5,10$. Thus the auxiliary logarithms involve only $2$ and $5$. Complex Lambert-series arguments are unavoidable, but the algebraic recursion still closes. For $p=11$, the new values include $17$ and $26=2\cdot13$: now the formula for $\log11$ introduces $\log13$ and $\log17$, primes greater than $11$.

Let $P(p)$ be the set of primes dividing $\prod_{n=1}^{(p-1)/2}(1+n^2)$. The following proposition states the observed boundary precisely.

\begin{proposition}[Closure and its verified range]\label{prop:closure}
The symmetric construction closes at $p=3,5,7$. It fails to close at every prime $p$ with $11\le p<20000$, in the sense that $P(p)$ contains a prime strictly larger than $p$.
\end{proposition}

For the first cases, this follows from elementary factorization. The remaining range is a computer-assisted finite verification by exact integer factorization; the complete reproducible procedure is given in Appendix~A. The proposition does not assert nonclosure for all primes above $7$.

\begin{table}[htbp]
\centering
\small
\setlength{\tabcolsep}{3.5pt}
\renewcommand{\arraystretch}{1.12}

\begin{tabular}{@{}rlccc@{}}
\toprule
$p$ & Auxiliary primes & Complex
arguments & Maximum $|q|$ & Closure \\
\midrule
\textbf{2} & ---                   & $0$  & $e^{-2\pi}$   & Start \\
\textbf{3} & $2$                   & $0$  & $e^{-2\pi}$   & Yes \\
\textbf{5} & $2$                   & $0$  & $e^{-2\pi}$   & Yes \\
7  & $2,5$                 & $4$  & $e^{-7\pi/5}$ & Yes \\
11 & $2,5,13,17$           & $8$  & --- & No \\
13 & $2,5,17,37$           & $8$  & --- & No \\
17 & $2,5,13,37$           & $12$ & --- & No \\
19 & $2,5,13,17,37,41$     & $16$ & --- & No \\
\bottomrule
\end{tabular}

\caption{Comparison of prime-logarithm identities
in the selected construction. For $2$, $3$, and $5$,
all Lambert-series arguments are real, with
$|q|\le e^{-2\pi}$, ensuring rapid convergence.
Prime $7$ marks the boundary: complex arguments
appear and the largest $|q|$ increases to
$e^{-7\pi/5}$, slowing convergence, but the
recursion still closes. For the displayed primes
$p\ge11$, closure fails. Auxiliary primes exclude
$p$ itself, and dashes indicate that no reduced
identity is considered. The identity for $\log2$
serves as the starting identity.}

\label{tab:closure}
\end{table}

Table~\ref{tab:closure} brings the primes $2$, $3$, and $5$
together: their selected identities use only real Lambert-series
arguments, and none has an argument larger in magnitude than
$e^{-2\pi}$. At $7$, complex arguments and a more slowly
converging series appear; as shown below, its leading term also
requires a much finer octave grid. The recursion still closes
at $7$, but fails for each tested prime beginning with $11$.
The general identity remains exact even when closure fails.

\subsection{An explicit illustration at prime $7$}

Eliminating the auxiliary $\log2$ and $\log5$ terms in \eqref{eq:general} gives the exact identity
\begin{align}\label{eq:log7}
\log7={}&\frac{781\pi}{1260}-\frac{80}{21}\mathcal L(e^{-2\pi})
+\frac{20}{21}\mathcal L(e^{-4\pi})
+\frac27\mathcal L(-e^{-5\pi})
+\frac47\mathcal L(-e^{-7\pi})\nonumber\\
&+\frac27\mathcal L(e^{-10\pi})+\frac47\mathcal L(e^{-14\pi})
+\frac47\operatorname{Re}\mathcal L(e^{-14\pi/5+2\pi i/5})\nonumber\\
&+\frac47\operatorname{Re}\mathcal L(e^{-7\pi/5+\pi i/5}).
\end{align}
The leading $\pi$-term is $781\pi/1260$, with the Lambert-series
terms providing the exact correction. Unlike the selected identities
for $\log2$, $\log3$, and $\log5$, whose series arguments have
absolute values no greater than $e^{-2\pi}$, the identity for
$\log7$ contains an argument of magnitude $e^{-7\pi/5}$.
This is approximately $6.6$ times larger than $e^{-2\pi}$,
so the corresponding Lambert series converges more slowly.
The convergence rate is therefore a further distinction of
these selected identities.

Relative to the leading term $2\pi/9$ for $\log2$, the
leading term for $\log7$ gives an approximate octave
coordinate of $781/280$. Incorporating this coordinate into
the $48$-step grid generated by $2$, $3$, and $5$ requires
\[
\operatorname{lcm}(48,280)=1680
\]
equal steps per octave, an increase by a factor of $35$.
This enlargement is a property of the present construction,
not a restriction on musical intervals involving $7$.

\subsection{What the arithmetic does, and does not, say about music}

The prime $7$ therefore occupies an intermediate position
in our construction. Compared with $2$, $3$, and $5$, it
requires complex Lambert-series arguments, includes a
more slowly converging series, and substantially enlarges
the common octave grid. Nevertheless, its identity still
closes using previously obtained prime logarithms. The
first three primes thus enjoy three forms of simplicity
that $7$ does not, while $7$ retains a property that fails
for the larger primes examined here.

This distinction invites comparison with musical traditions
in which $5$-limit tuning models are prominent. However,
the mathematics does not explain hearing, consonance, or
the historical development of tuning systems. Intervals
involving $7$ and higher primes remain musically legitimate.
The parallel is therefore a reason for further investigation,
not an explanation of musical practice.

\section{From prime logarithms to musical intervals}

\subsection{Why logarithms become pitches}

The logarithm is natural here because multiplying frequency ratios adds musical intervals. The cent value of a ratio $R$ is
\begin{equation}\label{eq:cents}
 C(R)=1200\log_2R=1200\frac{\log R}{\log2}.
\end{equation}
A pure octave has $1200$ cents. A $12$-tone equal-tempered semitone has $100$ cents. Thus any identity for $\log2,\log3,\log5$ immediately yields an identity for the pitch of every $5$-limit ratio in \eqref{eq:five_limit}.

\subsection{The $48$-step grid from the construction}

We now return to the three exact identities displayed together in
Section~2, equations~\eqref{eq:log2}--\eqref{eq:log5}. Their leading
$\pi$-terms provide the simple octave-grid approximation;
the remaining Lambert-series terms give the exact correction.

If we temporarily omit the Lambert-series remainders, the leading $\pi$-terms give
\begin{equation}\label{eq:pi_leading}
\widetilde{\log}\,2=\frac{2\pi}{9},\quad
\widetilde{\log}\,3=\frac{19\pi}{54},\quad
\widetilde{\log}\,5=\frac{37\pi}{72}.
\end{equation}
Consequently, the approximate fifth and major third are
\begin{equation}\label{eq:approxpitches}
\widetilde{\log}_2(3/2)=\frac7{12},\qquad
\widetilde{\log}_2(5/4)=\frac5{16}.
\end{equation}
These correspond to $700$ and $375$ cents, respectively. The fifth coincides with its position in $12$-tone equal temperament, while the major third lies $25$ cents below the equal-tempered value of $400$ cents. Since the least common multiple of $12$ and $16$ is $48$, both intervals belong to a common lattice of $48$ equal divisions of the octave.

The leading $\pi$-terms determine the assignments of the
three prime generators exactly. Expressing each in units of
one forty-eighth of the leading $\pi$-term for $\log2$ gives
\begin{equation}\label{eq:exactval}
\begin{aligned}
48\frac{2\pi/9}{2\pi/9}&=48,\\
48\frac{19\pi/54}{2\pi/9}&=76,\\
48\frac{37\pi/72}{2\pi/9}&=111.
\end{aligned}
\end{equation}
Thus the triplet $(48,76,111)$ follows directly from the
rational coefficients of $\pi$, without rounding.

This is also the familiar mapping of the prime generators in
$48$ equal divisions of the octave ($48$-EDO), a system with
$25$-cent steps that contains all $12$-TET pitches.\footnote{For
the $48$-EDO tuning system and its conventional prime mapping,
see the specialist \emph{Xenharmonic Wiki},
\url{https://en.xen.wiki/w/48edo}. This is a technical
reference for the tuning system, not for the Lambert-series
construction.} In conventional tuning theory, the same
integers arise by rounding the exact logarithmic coordinates:
\[
48\log_2 2=48,\qquad
48\log_2 3\approx76.0782,\qquad
48\log_2 5\approx111.4525.
\]
The Lambert-series corrections account exactly for the
differences, approximately $0.0782$ and $0.4525$ grid steps
for the generators $3$ and $5$. Both are smaller than half
a step, so the exact leading-term assignments agree with
nearest-step rounding. The agreement is especially close to
the rounding boundary for $5$, whose correction is about
$11.31$ cents compared with a half-step of $12.5$ cents.
This matching concerns the prime generators; their additive
combinations need not be the nearest grid points for every
compound interval, as illustrated below.

The Lambert-series identities therefore give a separate,
exact algebraic route to an established tuning map, while
their remainders recover the corresponding just intervals
and commas. The $22$ ratios of the shruti model in
\cite{Thakur2015} provide a further illustration
(Table~\ref{tab:shrutis}).

For a general $5$-limit interval, the resulting step count is
\begin{equation}\label{eq:grid}
 \widetilde C(2^a3^b5^c)=25(48a+76b+111c)
 \quad\text{cents}.
\end{equation}
The exponent of $2$ shifts the pitch by octaves; the exponent of $3$ contributes $76$ steps, equivalent to $28$ steps within one octave; the exponent of $5$ contributes $111$ steps, equivalent to $15$. This simple addition rule makes the approximation easy to use without evaluating any infinite sum.

In the next statement, $\mathbb Z$ denotes the integers and $\mathbb Z/48\mathbb Z$ denotes the $48$ possible grid positions, numbered $0$ through $47$, with positions differing by a whole octave identified.

\begin{proposition}[Pitch-class map]\label{prop:pitchmap}
The rule
\[
 \Phi:\mathbb Z^3\longrightarrow\mathbb Z/48\mathbb Z,\qquad
 \Phi(a,b,c)=28b+15c\pmod{48}
\]
is a surjective group homomorphism. Using only $2$ and $3$ gives exactly the twelve pitch classes divisible by $4$; allowing $5$ reaches all $48$ classes.
\end{proposition}

The musical interpretation is direct: if we ignore octaves, multiplying two frequency ratios adds their step numbers modulo $48$. Appendix~A supplies the short proof.

\subsection{The exact correction to a tempered note}

The $48$-step prediction is a convenient starting point, not the true pitch of a just interval. Write
\[
 \delta_3=\log3-\frac{19}{12}\log2,\qquad
 \delta_5=\log5-\frac{37}{16}\log2.
\]
These two quantities measure precisely how much the exact third- and fifth-prime logarithms differ from their grid assignments. Once they are known, \emph{every} $5$-limit interval is corrected by the same rule.

\begin{theorem}[Exact pitch decomposition]\label{thm:decomp}
For all integers $a,b,c$,
\begin{equation}\label{eq:decomp}
 C(2^a3^b5^c)=25(48a+76b+111c)
 +\frac{1200}{\log2}(b\delta_3+c\delta_5).
\end{equation}
The corrections have the exact Lambert-series forms
\begin{align}\label{eq:delta3}
\delta_3={}&\frac23\mathcal L(e^{-2\pi})
+\frac43\mathcal L(-e^{-3\pi})
-\frac{10}{3}\mathcal L(e^{-4\pi})
+\frac43\mathcal L(e^{-6\pi}),\\
\delta_5={}&\frac72\mathcal L(e^{-2\pi})
-\frac{11}{2}\mathcal L(e^{-4\pi})
+\mathcal L(-e^{-5\pi})
+\mathcal L(e^{-10\pi}).\label{eq:delta5}
\end{align}
\end{theorem}

Numerically, $1200\delta_3/\log2\approx1.955$ cents and $1200\delta_5/\log2\approx11.314$ cents. Consequently, $3/2$ lies about $1.955$ cents above the $700$-cent grid fifth, while $5/4$ lies about $11.314$ cents above the $375$-cent grid third. Notice the musical advantage of \eqref{eq:decomp}: the reader need only know the exponents $b$ and $c$. The proof is a substitution of the logarithmic identities and is given in Appendix~A.

\subsection{Two familiar commas as exact remainders}

The \emph{Pythagorean comma} measures the failure of twelve pure fifths to equal seven octaves. Since the $48$-step predictions close exactly, the entire comma is a Lambert-series remainder:
\begin{equation}\label{eq:pythcomma}
\begin{aligned}
\log\!\left(\frac{3^{12}}{2^{19}}\right)
={}&12\delta_3\\
={}&8\mathcal L(e^{-2\pi})+16\mathcal L(-e^{-3\pi})
-40\mathcal L(e^{-4\pi})+16\mathcal L(e^{-6\pi}).
\end{aligned}
\end{equation}
In pitch units this is approximately $23.46$ cents, or twelve increments of roughly $1.955$ cents.

The \emph{syntonic comma} compares the Pythagorean third $81/64$ with the pure third $5/4$. Their ratio is $81/80$; the selected lattice assigns one $25$-cent step to it, and the exact correction is negative:
\begin{equation}\label{eq:syntonic}
 \log\!\left(\frac{81}{80}\right)=\frac{\log2}{48}+4\delta_3-\delta_5.
\end{equation}
The correction is approximately $-3.49$ cents, leaving the true comma at about $21.51$ cents. To make the ``exact'' claim fully explicit, the remainder itself is
\begin{equation}\label{eq:syntonictail}
\begin{aligned}
4\delta_3-\delta_5={}&-\frac56\mathcal L(e^{-2\pi})
+\frac{16}{3}\mathcal L(-e^{-3\pi})
-\frac{47}{6}\mathcal L(e^{-4\pi})\\
&-\mathcal L(-e^{-5\pi})
+\frac{16}{3}\mathcal L(e^{-6\pi})
-\mathcal L(e^{-10\pi}).
\end{aligned}
\end{equation}
Both identities follow immediately from \eqref{eq:decomp}, as checked in Appendix~A. Thus the same two corrections account for discrepancies that arise in different parts of tuning theory.

\subsection{A scale example}

Table~\ref{tab:major_scale} shows a just major scale using the ratios discussed in connection with the \emph{Bilawal} model \cite{Thakur2015}, together with the tritone $45/32$. These illustrative ratios are not intended to prescribe fixed pitches in Hindustani performance. The table separates three quantities: the familiar $12$-TET position, our selected $48$-step approximation, and the exact just pitch.

\begin{table}[htbp]
\centering
\caption{A just major scale compared with $12$-TET
and the selected $48$-step mapping.
Here $\epsilon_{12}$ and $\epsilon_{48}$ are
the errors (just pitch minus assigned pitch),
in cents. Boldface identifies the more
accurate approximation when the two differ.}
\label{tab:major_scale}

\begingroup
\small
\setlength{\tabcolsep}{3pt}
\renewcommand{\arraystretch}{1.08}

\begin{tabular}{@{}lccrrrrr@{}}
\toprule
Note & Ratio & $(b,c)$ & $12$-TET &
$48$-EDO & Just &
$\epsilon_{12}$ & $\epsilon_{48}$ \\
\midrule
Do (Sa) & $1$ & $(0,0)$ &
0 & 0 & 0.00 & 0.00 & 0.00 \\

Re & $9/8$ & $(2,0)$ &
200 & 200 & 203.91 & $+3.91$ & $+3.91$ \\

Mi (Ga) & $5/4$ & $(0,1)$ &
400 & 375 & 386.31 & $-13.69$ &
$\mathbf{+11.31}$ \\

Fa (Ma) & $4/3$ & $(-1,0)$ &
500 & 500 & 498.04 & $-1.96$ & $-1.96$ \\

$\sharp$Fa & $45/32$ & $(2,1)$ &
600 & 575 & 590.22 &
$\mathbf{-9.78}$ & $+15.22$ \\

Sol (Pa) & $3/2$ & $(1,0)$ &
700 & 700 & 701.96 & $+1.96$ & $+1.96$ \\

La (Dha) & $27/16$ & $(3,0)$ &
900 & 900 & 905.87 & $+5.87$ & $+5.87$ \\

Ti (Ni) & $15/8$ & $(1,1)$ &
1100 & 1075 & 1088.27 &
$\mathbf{-11.73}$ & $+13.27$ \\

Do' & $2$ & $(0,0)$ &
1200 & 1200 & 1200.00 & 0.00 & 0.00 \\

\bottomrule
\end{tabular}
\endgroup
\end{table}

Table~\ref{tab:major_scale} shows that the finer $48$-step
grid does not always yield a closer approximation under the
selected additive assignment. It improves the just major
third $5/4$, while $12$-TET more closely approximates the
tritone $45/32$ and major seventh $15/8$. Several other
intervals are represented equally well. The familiar
$12$-tone system thus provides effective approximations
despite its coarser resolution. This does not mean the full
$48$-EDO grid is less accurate: it contains every $12$-TET
position. The difference arises because our assignment
adds the contributions of the prime generators rather than
rounding every interval independently.

\subsection{The $22$ shrutis on the grid}

The shruti model of \cite{Thakur2015} assigns $22$ ratios to the octave. Each has the form $2^a3^b5^c$ with $c\in\{0,\pm1\}$, so \eqref{eq:grid} places it on the $48$-step lattice directly. Table~\ref{tab:shrutis} lists all $22$. They occupy $22$ of the $48$ grid steps, each within $15.3$ cents of its grid point. They split as $12$ ratios with $c=0$ (the Pythagorean chain), $5$ with $c=+1$, and $5$ with $c=-1$.

\begin{table}[htbp]
\centering
\caption{The $22$ shruti ratios of the model
in \cite{Thakur2015}, mapped onto the
$48$-step grid generated by the primes
$2$, $3$, and $5$. All $22$ occupy distinct
grid positions. Here $(b,c)$ are the exponents
of $3$ and $5$, $k=28b+15c\pmod{48}$ gives
the grid position, and
$\Delta C=b\Delta_3+c\Delta_5$ gives the
signed correction in cents, with
$\Delta_3\approx1.955$ and
$\Delta_5\approx11.314$.
Shruti labels follow \cite{Thakur2015}.}

\label{tab:shrutis}

\begingroup
\small
\setlength{\tabcolsep}{5pt}
\renewcommand{\arraystretch}{1.08}

\begin{tabular}{@{}lccrr@{}}
\toprule
Shruti & Ratio & $(b,c)$ & Step $k$ &
$\Delta C$ (cents) \\
\midrule
S         & $1$       & $(0,0)$   & 0  & $0.00$ \\
R$_{11}$  & $256/243$ & $(-5,0)$  & 4  & $-9.78$ \\
R$_{12}$  & $16/15$   & $(-1,-1)$ & 5  & $-13.27$ \\
R$_{21}$  & $10/9$    & $(-2,1)$  & 7  & $+7.40$ \\
R$_{22}$  & $9/8$     & $(2,0)$   & 8  & $+3.91$ \\
G$_{11}$  & $32/27$   & $(-3,0)$  & 12 & $-5.87$ \\
G$_{12}$  & $6/5$     & $(1,-1)$  & 13 & $-9.36$ \\
G$_{21}$  & $5/4$     & $(0,1)$   & 15 & $+11.31$ \\
G$_{22}$  & $81/64$   & $(4,0)$   & 16 & $+7.82$ \\
M$_{11}$  & $4/3$     & $(-1,0)$  & 20 & $-1.96$ \\
M$_{12}$  & $27/20$   & $(3,-1)$  & 21 & $-5.45$ \\
M$_{21}$  & $45/32$   & $(2,1)$   & 23 & $+15.22$ \\
M$_{22}$  & $729/512$ & $(6,0)$   & 24 & $+11.73$ \\
P         & $3/2$     & $(1,0)$   & 28 & $+1.96$ \\
D$_{11}$  & $128/81$  & $(-4,0)$  & 32 & $-7.82$ \\
D$_{12}$  & $8/5$     & $(0,-1)$  & 33 & $-11.31$ \\
D$_{21}$  & $5/3$     & $(-1,1)$  & 35 & $+9.36$ \\
D$_{22}$  & $27/16$   & $(3,0)$   & 36 & $+5.87$ \\
N$_{11}$  & $16/9$    & $(-2,0)$  & 40 & $-3.91$ \\
N$_{12}$  & $9/5$     & $(2,-1)$  & 41 & $-7.40$ \\
N$_{21}$  & $15/8$    & $(1,1)$   & 43 & $+13.27$ \\
N$_{22}$  & $243/128$ & $(5,0)$   & 44 & $+9.78$ \\
\bottomrule
\end{tabular}
\endgroup
\end{table}
\FloatBarrier

The $48$-step mapping preserves distinctions that are
lost in $12$-tone equal temperament. Rounding the $22$
theoretical shruti ratios to their nearest $12$-TET
notes gives only $12$ distinct pitches, whereas the
additive $48$-step mapping keeps all $22$ distinct.
Each of the ten shrutis with $c=\pm1$ differs by
exactly one syntonic comma $81/80$ from a Pythagorean
shruti with $c=0$. The additive mapping places all
ten pairs on adjacent steps. By contrast, independent
rounding to the nearest $25$ cents merges three pairs:
$256/243$ and $16/15$, $45/32$ and $729/512$,
and $15/8$ and $243/128$.

The mean absolute deviation from the exact ratios is
approximately $9.86$ cents for nearest-note $12$-TET
and $7.90$ cents for the present $48$-step mapping.
The largest deviations are $19.55$ and $15.22$ cents,
respectively. The $48$-step assignment is closer for
seven ratios, $12$-TET for three, and both are equally
accurate for the remaining twelve. Thus the finer grid
can preserve distinctions between nearby shrutis that
are lost in the smaller system, although $12$-TET
remains an effective coarser approximation.

These results illustrate the usefulness of the $48$-step
construction for representing a theoretical system of
closely spaced pitches. They do not establish that this
grid is musically preferable to $12$-TET, nor do they
prescribe the intonation of actual performances, which
varies with musical context
\cite{JairazbhoyStone1963,SerraEtAl2011}.

\subsection{An octave-circle picture}

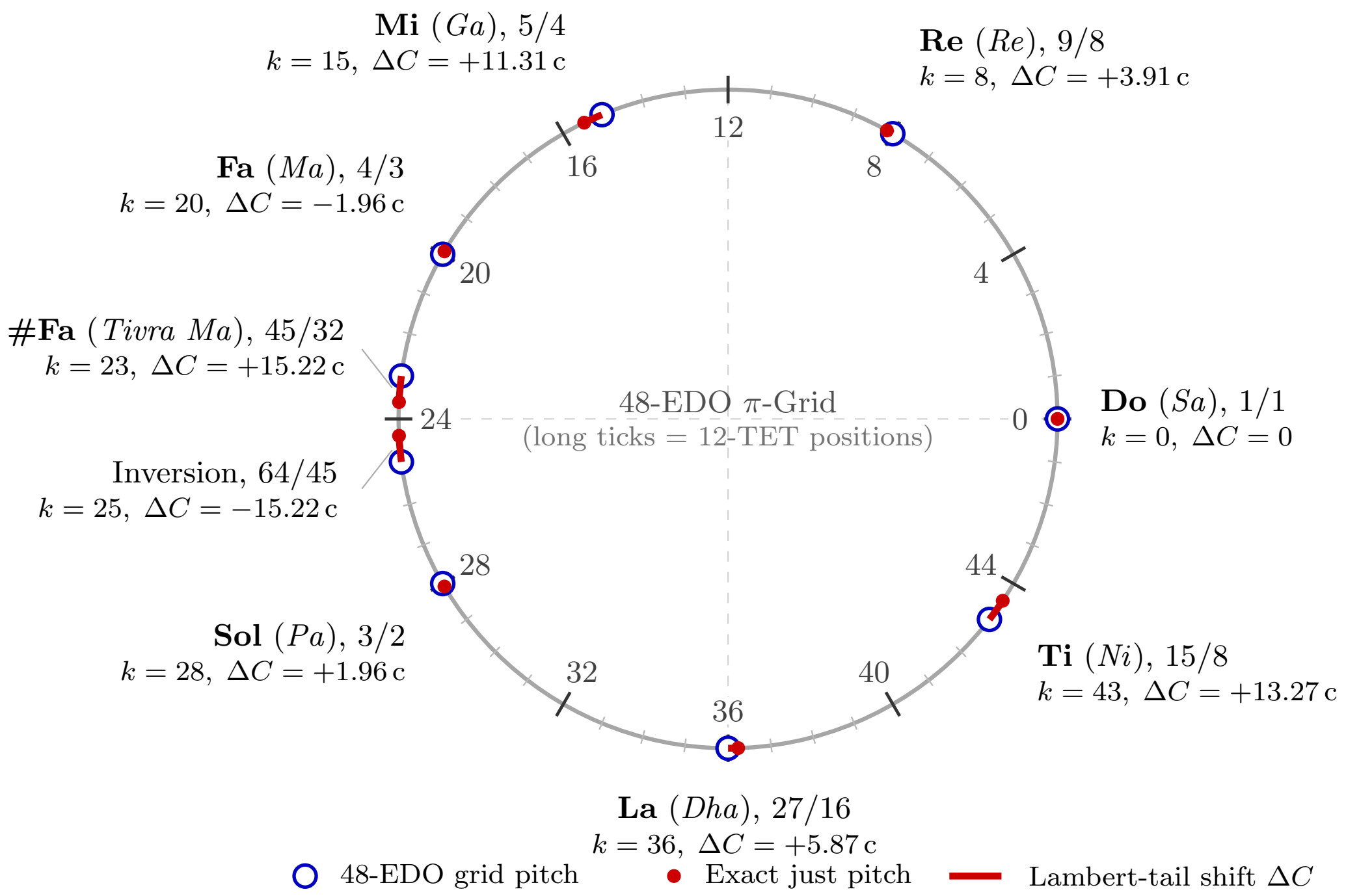
\begin{figure}[H]
\centering
\begin{tikzpicture}
    \def\R{3.35}

    \draw[gray!35, dashed] (-\R,0) -- (\R,0);
    \draw[gray!35, dashed] (0,-\R) -- (0,\R);

    \draw[gray!70, line width=1.2pt] (0,0) circle (\R);

    \foreach \k in {0,1,...,47} {
        \draw[gray!55, line width=0.45pt]
            ({\k*7.5}:{\R-0.07}) -- ({\k*7.5}:{\R+0.07});
    }

    \foreach \k in {0,4,...,44} {
        \draw[black!80, line width=0.9pt]
            ({\k*7.5}:{\R-0.14}) -- ({\k*7.5}:{\R+0.14});
        \node[font=\footnotesize, text=black!75, fill=white, inner sep=1pt]
            at ({\k*7.5}:{\R-0.38}) {$\k$};
    }

    \node[align=center, font=\footnotesize, text=black!70] at (0,0.18)
        {$48$-EDO $\pi$-Grid};
    \node[align=center, font=\scriptsize, text=black!55] at (0,-0.20)
        {(long ticks = $12$-TET positions)};

    \draw[blue!75!black, fill=white, line width=1.0pt] (0:\R) circle (3.2pt);
    \fill[red!80!black] (0:\R) circle (2.0pt);
    \node[anchor=west, align=left, font=\footnotesize] at (0:{\R+0.30})
        {\textbf{Do} (\emph{Sa}), $1/1$\\[-1pt]
        {\scriptsize $k=0,\;\Delta C=0$}};

    \draw[blue!75!black, fill=white, line width=1.0pt] (60:\R) circle (3.2pt);
    \draw[red!80!black, line width=2.0pt] (60:\R) arc (60:61.17:\R);
    \fill[red!80!black] (61.17:\R) circle (2.0pt);
    \node[anchor=south west, align=left, font=\footnotesize] at (60.6:{\R+0.33})
        {\textbf{Re} (\emph{Re}), $9/8$\\[-1pt]
        {\scriptsize $k=8,\;\Delta C=+3.91$\,c}};

    \draw[blue!75!black, fill=white, line width=1.0pt] (112.5:\R) circle (3.2pt);
    \draw[red!80!black, line width=2.0pt] (112.5:\R) arc (112.5:115.89:\R);
    \fill[red!80!black] (115.89:\R) circle (2.0pt);
    \node[anchor=south east, align=right, font=\footnotesize] at (114.2:{\R+0.34})
        {\textbf{Mi} (\emph{Ga}), $5/4$\\[-1pt]
        {\scriptsize $k=15,\;\Delta C=+11.31$\,c}};

    \draw[blue!75!black, fill=white, line width=1.0pt] (150:\R) circle (3.2pt);
    \draw[red!80!black, line width=2.0pt] (150:\R) arc (150:149.41:\R);
    \fill[red!80!black] (149.41:\R) circle (2.0pt);
    \node[anchor=south east, align=right, font=\footnotesize] at (148.7:{\R+0.34})
        {\textbf{Fa} (\emph{Ma}), $4/3$\\[-1pt]
        {\scriptsize $k=20,\;\Delta C=-1.96$\,c}};

    \draw[blue!75!black, fill=white, line width=1.0pt] (172.5:\R) circle (3.2pt);
    \draw[red!80!black, line width=2.0pt] (172.5:\R) arc (172.5:177.07:\R);
    \fill[red!80!black] (177.07:\R) circle (2.0pt);
    \draw[gray!65, thin] (174.7:{\R+0.08}) -- (169.2:{\R+0.44});
    \node[anchor=east, align=right, font=\footnotesize] at (169.2:{\R+0.48})
        {\textbf{\#Fa} (\emph{Tivra Ma}), $45/32$\\[-1pt]
        {\scriptsize $k=23,\;\Delta C=+15.22$\,c}};

    \draw[blue!75!black, fill=white, line width=1.0pt] (187.5:\R) circle (3.2pt);
    \draw[red!80!black, line width=2.0pt] (187.5:\R) arc (187.5:182.93:\R);
    \fill[red!80!black] (182.93:\R) circle (2.0pt);
    \draw[gray!65, thin] (185.3:{\R+0.08}) -- (190.8:{\R+0.44});
    \node[anchor=east, align=right, font=\footnotesize] at (190.8:{\R+0.55})
        {Inversion, $64/45$\\[-1pt]
        {\scriptsize $k=25,\;\Delta C=-15.22$\,c}};

    \draw[blue!75!black, fill=white, line width=1.0pt] (210:\R) circle (3.2pt);
    \draw[red!80!black, line width=2.0pt] (210:\R) arc (210:210.59:\R);
    \fill[red!80!black] (210.59:\R) circle (2.0pt);
    \node[anchor=north east, align=right, font=\footnotesize] at (211.6:{\R+0.34})
        {\textbf{Sol} (\emph{Pa}), $3/2$\\[-1pt]
        {\scriptsize $k=28,\;\Delta C=+1.96$\,c}};

    \draw[blue!75!black, fill=white, line width=1.0pt] (270:\R) circle (3.2pt);
    \draw[red!80!black, line width=2.0pt] (270:\R) arc (270:271.76:\R);
    \fill[red!80!black] (271.76:\R) circle (2.0pt);
    \node[anchor=north, align=center, font=\footnotesize] at (271:{\R+0.34})
        {\textbf{La} (\emph{Dha}), $27/16$\\[-1pt]
        {\scriptsize $k=36,\;\Delta C=+5.87$\,c}};

    \draw[blue!75!black, fill=white, line width=1.0pt] (322.5:\R) circle (3.2pt);
    \draw[red!80!black, line width=2.0pt] (322.5:\R) arc (322.5:326.48:\R);
    \fill[red!80!black] (326.48:\R) circle (2.0pt);
    \node[anchor=north west, align=left, font=\footnotesize] at (324.6:{\R+0.34})
        {\textbf{Ti} (\emph{Ni}), $15/8$\\[-1pt]
        {\scriptsize $k=43,\;\Delta C=+13.27$\,c}};

    \begin{scope}[shift={(0,-\R-1.30)}]
        \draw[blue!75!black, fill=white, line width=1.0pt] (-4.3,0) circle (3.2pt);
        \node[anchor=west, font=\scriptsize] at (-4.08,0)
            {$48$-EDO grid pitch};

        \fill[red!80!black] (-0.55,0) circle (2.0pt);
        \node[anchor=west, font=\scriptsize] at (-0.37,0)
            {Exact just pitch};

        \draw[red!80!black, line width=2.0pt] (2.25,0) -- (2.78,0);
        \node[anchor=west, font=\scriptsize] at (2.92,0)
            {Lambert-tail shift $\Delta C$};
    \end{scope}
\end{tikzpicture}
\caption{An octave circle showing the selected $48$-EDO grid
(open blue circles), exact just pitches (red dots), and
Lambert-series corrections (red arcs). Longer ticks mark the
$12$-TET positions. Labels show Western and Indian note names,
just ratios, grid indices $k$, and the signed pitch corrections
$\Delta C$. The tritone and its octave inversion have equal and
opposite corrections.}
\label{fig:pitch_circle}
\end{figure}

Figure~\ref{fig:pitch_circle} shows the same information geometrically. Blue open circles locate the selected grid pitches, red points locate exact just pitches, and short arcs visualize the signed corrections. Inversion of an interval about the octave reverses the correction's sign.

The figure illustrates mathematical pitch relationships rather than prescribing musical performance. Perceived consonance also depends on sound spectra and hearing \cite{PlompLevelt1965,McDermottOxenham2008}.
\FloatBarrier

\section{Discussion and conclusions}

The construction distinguishes $2$, $3$, and $5$ in three concrete ways. Their selected prime-logarithm identities use only real Lambert-series arguments, all with magnitudes at most $e^{-2\pi}$, and their leading terms generate a compact $48$-step octave grid. The first two properties simplify the identities and favor rapid convergence; the third connects their leading logarithmic terms to an established tuning lattice. Prime $7$ marks a clear change: complex arguments appear, one series decays more slowly, and including its leading coordinate enlarges the grid to $1680$ steps. Yet its identity still closes using the earlier prime logarithms. For every prime from $11$ through $19997$, the finite exact computation finds an auxiliary prime larger than the prime being treated, so the same recursive closure fails.

The musical consequences are equally explicit. The $48$-step assignments give simple approximations to $5$-limit intervals, while the Lambert-series remainders recover their exact pitches and familiar commas. The additive $48$-step mapping distinguishes all $22$ ratios in the theoretical shruti model studied here, including pairs a syntonic comma apart; the comparison with $12$-tone equal temperament also shows that finer resolution does not make every individual additive assignment more accurate. These examples demonstrate both the usefulness and the limits of the construction as a framework for comparing tuning systems.

These results establish precise arithmetic distinctions that parallel the prominence of $5$-limit tuning. Whether this parallel has any perceptual or historical significance is a separate question requiring musical evidence. A natural direction for further work is a systematic study of how this grid organizes theoretical shrutis, scales, and ragas.

\appendix
\section{Proofs and computational details}

\subsection{Derivation of the general prime identity}

\begin{proof}[Proof of Theorem~\ref{thm:prime}]
For $n=-r,\ldots,r$, put $t_n=(1+in)/p$. The numbers $n$ represent every residue class modulo $p$, so prime splitting at $q=e^{-2\pi}$ gives
\begin{equation}\label{eq:appsplit}
\sum_{n=-r}^r\mathcal L(e^{-2\pi t_n})
=(p+1)\mathcal L(e^{-2\pi})-\mathcal L(e^{-2\pi p}).
\end{equation}
By modular inversion \eqref{eq:inversion},
\[
\mathcal L(e^{-2\pi t_n})
=\mathcal L(e^{-2\pi/t_n})
+\frac12\log t_n-\frac{\pi}{12}(t_n-t_n^{-1}).
\]
Pairing $n$ with $-n$ using the principal logarithm gives
\[
\sum_{n=-r}^r\log t_n
=\sum_{n=1}^r\log(1+n^2)-p\log p.
\]
The other elementary sums are
\[
\sum_{n=-r}^r t_n=1,\qquad
\sum_{n=-r}^r\frac1{t_n}=pS_p.
\]
Finally, $1/t_n=p/(1+in)$, so $e^{-2\pi/t_n}=z_{p,n}$. Summing the transformed series therefore yields
\[
\sum_{n=-r}^r\mathcal L(z_{p,n}).
\]
Substitute these three expressions into \eqref{eq:appsplit} and solve for $\log p$. This is exactly \eqref{eq:general}.
\end{proof}

\subsection{Real-argument criterion}

\begin{proof}[Proof of Theorem~\ref{thm:real}]
For any integer $n$,
\[
 z_{p,n}=\exp\!\left(-\frac{2\pi p}{1+n^2}\right)
 \exp\!\left(\frac{2\pi i pn}{1+n^2}\right).
\]
This is real exactly when $2pn/(1+n^2)$ is an integer. Since $\gcd(n,1+n^2)=1$, the condition is equivalent to $(1+n^2)\mid2p$.

For $p=2$, representatives $0,1$ suffice (this is the transformation at $t=1/2$ and $t=(1\pm i)/2$ that underlies \eqref{eq:log2}); for $p=3$, choose $-1,0,1$; and for $p=5$, choose $-2,-1,0,1,2$. In each case all $1+n^2$ divide $2p$.

If $p\ge7$, every complete residue system contains an integer $n\equiv2\pmod p$. Then $|n|\ge2$, so $1+n^2\ge5$, while $1+n^2\equiv5\pmod p$ and hence $p\nmid(1+n^2)$. The only positive divisors of $2p$ are $1,2,p,2p$, none of which can equal $1+n^2$. Thus one argument must be non-real.
\end{proof}

\subsection{Proof of the pitch-class map}

\begin{proof}[Proof of Proposition~\ref{prop:pitchmap}]
Reducing $48a+76b+111c$ modulo $48$ gives $28b+15c$. Multiplication of ratios adds exponents, hence adds classes. Since $28\cdot7+15\cdot(-13)=1$, the image contains $1$ and is all of $\mathbb Z/48\mathbb Z$. With $c=0$, the image consists exactly of the multiples of $\gcd(28,48)=4$, giving twelve classes.
\end{proof}

\subsection{Proof of the exact pitch decomposition and commas}

\begin{proof}[Proof of Theorem~\ref{thm:decomp}]
Use $\log3=(19/12)\log2+\delta_3$ and $\log5=(37/16)\log2+\delta_5$ in the logarithm of $2^a3^b5^c$. The coefficient of $\log2$ is $(48a+76b+111c)/48$. Multiplying by $1200/\log2$ yields \eqref{eq:decomp}. Subtracting the selected identities \eqref{eq:log2}--\eqref{eq:log5} gives \eqref{eq:delta3} and \eqref{eq:delta5} term by term.
\end{proof}

For the Pythagorean comma, set $(a,b,c)=(-19,12,0)$ in \eqref{eq:decomp}; the grid contribution is zero, leaving $12\delta_3$. For the syntonic comma, set $(a,b,c)=(-4,4,-1)$; the grid contributes $25$ cents, or $\log2/48$ in logarithmic units, leaving $4\delta_3-\delta_5$. Expanding both expressions proves \eqref{eq:pythcomma}--\eqref{eq:syntonictail}.

\subsection{Finite verification of the closure proposition}

\begin{proof}[Proof of Proposition~\ref{prop:closure}]
For $p=3$, the only value $1+n^2$ with $1\le n\le(p-1)/2$ is $2$. For $p=5$ these values are $2,5$; for $p=7$ they are $2,5,10$. All auxiliary prime factors are smaller than $p$, apart from $p$ itself in the $p=5$ case. Its coefficient can be moved to the left, so the recursion closes.

For $p\ge11$ in the stated finite range, the following exact-integer calculation checks whether some $1+n^2$ with $n\le(p-1)/2$ has a prime factor greater than $p$. Since the upper bound on $n$ increases with $p$, we can retain the largest prime factor already found. The assertion is the complete computational certificate for the finite range; it does not constitute a proof for unbounded $p$.
\end{proof}

\begin{verbatim}
from sympy import primerange, factorint
largest = 0
n = 0
for p in primerange(11, 20000):
    r = (p - 1) // 2
    while n < r:
        n += 1
        largest = max(largest, max(factorint(n*n + 1)))
    assert largest > p
print("All primes in the stated range pass.")
\end{verbatim}

\section*{Acknowledgments}

The author acknowledges the use of generative AI tools---Google Gemini, OpenAI ChatGPT, and Anthropic Claude---for exploratory calculations, LaTeX and TikZ assistance, and editorial refinement. The author takes responsibility for the mathematical results, their verification, and the final manuscript.

\end{document}